\documentclass[letterpaper, 10 pt, conference]{ieeeconf}
\IEEEoverridecommandlockouts
\usepackage{graphics}
\usepackage{epsfig}
\usepackage{epstopdf}
\usepackage{stfloats}
\usepackage{amsmath}
\usepackage{subcaption}
\usepackage{amssymb}
\usepackage{gensymb}
\usepackage{cite}
\usepackage{float}
\usepackage{hhline}
\usepackage{graphicx}
\usepackage{mathtools}
\usepackage{array,tabularx}
\usepackage{cellspace}
\usepackage{color}
\usepackage{algpseudocode}
\usepackage{xpatch}
\usepackage{tabu}
\usepackage{lipsum}
\usepackage[normalem]{ulem}
\newtheorem{problem}{Problem}
\newtheorem{algorithm}{Algorithm}

\newcommand{\boxing}[5]{
	\begin{figure}[#1]
	    \vspace{#3}
		\begin{center}
			\fbox{
				\parbox{#2\textwidth}{
					\vspace{0.2cm}
					#5
					\vspace{0.2cm}
				}
			}
		\end{center}
		\vspace{#4}\null
	\end{figure}
}

\newcommand{\sref}[1]{\S\ref{#1}}
\newcommand{\blkdiag}[1]{\mathop{\bf blkdiag}\left\lbrace #1 \right\rbrace}

\newcommand{\mat}[1]{\begin{bmatrix} #1 \end{bmatrix}}
\newcommand{\transp}{{\scriptscriptstyle\mathsf{T}}}
\newcommand{\psd}[1]{\mathbb{S}^{#1}_{+}}
\newcommand{\real}{\mathbb{R}}
\newcommand{\set}[1]{\mathcal{#1}}
\newcommand{\definedas}{\coloneqq}
\newcommand{\ith}[2]{#1^\textit{#2}}

\newcommand{\Nz}{{n_z}}
\newcommand{\Nc}{{n_c}}
\newcommand{\Ng}{{n_g}}
\newcommand{\shat}{\hat{\sigma}}

\newcommand{\tf}{t_f}
\newcommand{\kd}{k_d}
\newcommand{\thetamax}{\theta_\textit{max}}

\newcommand{\vmax}{v_\textit{max}}
\newcommand{\Tmin}{T_\textit{min}}
\newcommand{\Tmax}{T_\textit{max}}

\newcommand{\rh}{r_{h}}
\newcommand{\nh}{\hat{n}_{h}}
\newcommand{\Lh}{\ell_{c}}
\newcommand{\Rh}{\rho_{h}}
\newcommand{\Rhg}{\rho_{g}}
\newcommand{\Rhc}{\rho_{c}}
\newcommand{\Nho}{\hat{N}_{h}}

\newcommand{\Lspacing}{\ell_{o}}
\newcommand{\No}{N_o}
\newcommand{\Ro}{R_o}
\newcommand{\obset}{\set{N}_{o}}
\newcommand{\phato}{\hat{p}_{o}}
\newcommand{\qhato}{\hat{q}_{o}}
\newcommand{\Wo}{w_{o}}
\newcommand{\jo}{{l}}
\newcommand{\roj}{r_{o,\jo}}

\newcommand{\Nx}{{n_x}}
\newcommand{\Nu}{{n_u}}
\newcommand{\Nw}{{n_w}}
\newcommand{\Na}{{n_\alpha}}
\newcommand{\KK}{K}
\newcommand{\Kset}{\set{K}}
\newcommand{\Ksetm}{\bar{\set{K}}}
\newcommand{\dt}{\Delta t}
\newcommand{\tk}{t_{k}}
\newcommand{\xdk}[1]{x_{#1}}
\newcommand{\udk}[1]{u_{#1}}
\newcommand{\adk}[1]{\alpha_{#1}}
\newcommand{\zdk}[1]{z_{#1}}
\newcommand{\vdk}[1]{\nu_{#1}}

\newcommand{\ud}{\bar{\udk{}}}
\newcommand{\zd}{\bar{\zdk{}}}
\newcommand{\vd}{\bar{\vdk{}}}
\newcommand{\Ad}{A_d}
\newcommand{\Bd}{B_d}
\newcommand{\Ed}{E_d}
\newcommand{\Wd}{w_d}
\newcommand{\Wvc}{W_{vc}}
\newcommand{\Wtr}{W_{tr}}
\newcommand{\Jorig}{J}
\newcommand{\Jvc}{J_\textit{vc}}
\newcommand{\Jtr}{J_\textit{tr}}
\newcommand{\vctol}{\epsilon_\textit{vc}}
\newcommand{\trtol}{\epsilon_\textit{tr}}

\usepackage{rotating}
\usepackage{pgfplots}
\usepackage{tikz}
\usetikzlibrary{patterns}
\usetikzlibrary{math}
\usetikzlibrary{scopes}
\usetikzlibrary{fadings}
\usetikzlibrary{arrows,calc,decorations.pathmorphing}
\usetikzlibrary{matrix,positioning}
\tikzset{>=latex}

\definecolor{beige}{RGB}{245,245,220}
\definecolor{darkred}{rgb}{0.90,0.00,0.00}
\definecolor{darkgreen}{rgb}{0.00,0.45,0.00}
\definecolor{midgreen}{rgb}{0.00,0.55,0.00}

\tikzset{ 
table/.style={
  matrix of math nodes,
  row sep=-\pgflinewidth,
  column sep=-\pgflinewidth,
  nodes={rectangle,draw=white,text width=6em,align=left},
  text depth=0.25ex,
  text height=2ex,
  nodes in empty cells
  },
  title/.style={font=\large}
}

\definecolor{ucol}{RGB}{255,0,0}
\definecolor{gcol}{RGB}{0,120,0}
\definecolor{scol}{RGB}{63, 226, 45}
\definecolor{vcol}{RGB}{0,0,0}
\definecolor{acol}{RGB}{0, 128, 255}
\definecolor{dcol}{RGB}{204,102,0}
\definecolor{lcol}{RGB}{204,102,0}
\definecolor{bcol}{RGB}{0,0,0}
\definecolor{ocol}{RGB}{167,167,167}

\newcommand{\threeaxes}[8]{
	\tikzmath{
    	\Lxr=#3;\Lxl=#3;\Lyt=#4;\Lyb=#4;\Zt= #5;\Zb= #6;\Xang=#7;\Yang=#8;
    	\Xxr= cos(\Xang)*\Lxr; \Xyr=-sin(\Xang)*\Lxr;
    	\Xxl=-cos(\Xang)*\Lxl; \Xyl= sin(\Xang)*\Lxl;
        \Yxt= sin(\Yang)*\Lyt; \Yyt= cos(\Yang)*\Lyt;
        \Yxb=-sin(\Yang)*\Lyb; \Yyb=-cos(\Yang)*\Lyb;
        \zzz=0;
    }
    \begin{scope}[shift={(#1,#2)},rotate=0]
        \ifx\Zt\zzz\else  \draw[black,->] (0,0) -- +(0, \Zt);    \fi
        \ifx\Zb\zzz\else  \draw[black,->] (0,0) -- +(0,-\Zb);    \fi
        \ifx\Lxr\zzz\else \draw[black,->] (0,0) -- +(\Xxr,\Xyr); \fi
        \ifx\Lxl\zzz\else \draw[black,->] (0,0) -- +(\Xxl,\Xyl); \fi
        \ifx\Lyt\zzz\else \draw[black,->] (0,0) -- +(\Yxt,\Yyt); \fi
        \ifx\Lyb\zzz\else \draw[black,->] (0,0) -- +(\Yxb,\Yyb); \fi
	\end{scope}
}
\newcommand{\threeaxeslabelx}[9]{
	\tikzmath{
    	\Lxr=#3;\Lxl=#3;\Lyt=#4;\Lyb=#4;\Zt= #5;\Zb= #6;\Xang=#7;\Yang=#8;
    	\Xxr= cos(\Xang)*\Lxr; \Xyr=-sin(\Xang)*\Lxr;
    	\Xxl=-cos(\Xang)*\Lxl; \Xyl= sin(\Xang)*\Lxl;
        \Yxt= sin(\Yang)*\Lyt; \Yyt= cos(\Yang)*\Lyt;
        \Yxb=-sin(\Yang)*\Lyb; \Yyb=-cos(\Yang)*\Lyb;
        \zzz=0;
    }
    \begin{scope}[shift={(#1,#2)},rotate=0]
        \ifx\Lxr\zzz\else \draw (\Xxr,\Xyr) node[anchor=west] {#9}; \fi
	\end{scope}
}
\newcommand{\threeaxeslabely}[9]{
	\tikzmath{
    	\Lxr=#3;\Lxl=#3;\Lyt=#4;\Lyb=#4;\Zt= #5;\Zb= #6;\Xang=#7;\Yang=#8;
    	\Xxr= cos(\Xang)*\Lxr; \Xyr=-sin(\Xang)*\Lxr;
    	\Xxl=-cos(\Xang)*\Lxl; \Xyl= sin(\Xang)*\Lxl;
        \Yxt= sin(\Yang)*\Lyt; \Yyt= cos(\Yang)*\Lyt;
        \Yxb=-sin(\Yang)*\Lyb; \Yyb=-cos(\Yang)*\Lyb;
        \zzz=0;
    }
    \begin{scope}[shift={(#1,#2)},rotate=0]
        \ifx\Lyt\zzz\else \draw (\Yxt,\Yyt) node[anchor=east] {#9}; \fi
	\end{scope}
}
\newcommand{\threeaxeslabelz}[9]{
	\tikzmath{
    	\Lxr=#3;\Lxl=#3;\Lyt=#4;\Lyb=#4;\Zt= #5;\Zb= #6;\Xang=#7;\Yang=#8;
    	\Xxr= cos(\Xang)*\Lxr; \Xyr=-sin(\Xang)*\Lxr;
    	\Xxl=-cos(\Xang)*\Lxl; \Xyl= sin(\Xang)*\Lxl;
        \Yxt= sin(\Yang)*\Lyt; \Yyt= cos(\Yang)*\Lyt;
        \Yxb=-sin(\Yang)*\Lyb; \Yyb=-cos(\Yang)*\Lyb;
        \zzz=0;
    }
    \begin{scope}[shift={(#1,#2)},rotate=0]
        \ifx\Zt\zzz\else  \draw (0,\Zt)     node[anchor=west] {$\;$#9}; \fi
	\end{scope}
}

\newcommand{\pane}[7]{
	\tikzmath{
    	\rot=#3;
    	\width=#4;
        \height=#5;
        \corner=#6;
    	\px1= 0.5*\width-\corner; \py1=-0.5*\height;
        \px2= 0.5*\width;         \py2=-0.5*\height+\corner;
        \px3= 0.5*\width;         \py3= 0.5*\height-\corner;
        \px4= 0.5*\width-\corner; \py4= 0.5*\height;
        \px5=-0.5*\width+\corner; \py5= 0.5*\height;
        \px6=-0.5*\width;         \py6= 0.5*\height-\corner;
        \px7=-0.5*\width;         \py7=-0.5*\height+\corner;
        \px8=-0.5*\width+\corner; \py8=-0.5*\height;
    }
	\begin{scope}[shift={(#1,#2)},rotate=\rot]
    	\filldraw[{#7}]
        (\px1,\py1) to[out=    0,in=  -90] (\px2,\py2) --
        (\px3,\py3) to[out=   90,in=    0] (\px4,\py4) --
        (\px5,\py5) to[out= -180,in=   90] (\px6,\py6) --
        (\px7,\py7) to[out=  -90,in= -180] (\px8,\py8) -- cycle;
    \end{scope}
}

\newcommand{\quadrotor}[6]{
	\tikzmath{
		\Larm = 1; \Rr = 0.9;
		\rx1= \Larm; \ry1= \Larm;
		\rx2= \Larm; \ry2=-\Larm;
		\rx3=-\Larm; \ry3=-\Larm;
		\rx4=-\Larm; \ry4= \Larm;
		\drx1=\rx1+\Rr*cos(45); \dry1=\ry1+\Rr*sin(45);
		\val3 = #3;
		\val2 = #4;
		\val1 = #5;
    \balllinewidth=0.20mm;
    \rotorlinewidth=0.40mm;
    \rotoroutercolor=40;
	}
	\begin{scope}[shift={({(#1)/#6},{(#2)/#6})}]
		\begin{scope}[transform canvas={scale=#6}]
			\begin{scope}[rotate=\val1]
				\begin{scope}[xscale=\val2]
					\begin{scope}[rotate=-\val3]
							\begin{scope}[shift={(\rx1,\ry1)}]
								\fill[shading=radial,outer color=orange!\rotoroutercolor,inner color=orange!100] (0,0) circle (\Rr);
								\draw[-,>=stealth,line cap=round,line width=\rotorlinewidth] ( 90:\Rr) arc ( 90:-270:\Rr);
							\end{scope}
							\begin{scope}[shift={(\rx2,\ry2)}]
								\fill[shading=radial,outer color=black!5,inner color=black!\rotoroutercolor] (0,0) circle (\Rr);
								\draw[-,>=stealth,line cap=round,line width=\rotorlinewidth] (-90:\Rr) arc (-90: 270:\Rr);
							\end{scope}
							\begin{scope}[shift={(\rx3,\ry3)}]
								\fill[shading=radial,outer color=black!5,inner color=black!\rotoroutercolor] (0,0) circle (\Rr);
								\draw[-,>=stealth,line cap=round,line width=\rotorlinewidth] (-90:\Rr) arc (-90: 270:\Rr);
							\end{scope}
							\begin{scope}[shift={(\rx4,\ry4)}]
								\fill[shading=radial,outer color=orange!\rotoroutercolor,inner color=orange!100] (0,0) circle (\Rr);
								\draw[-,>=stealth,line cap=round,line width=\rotorlinewidth] ( 90:\Rr) arc ( 90:-270:\Rr);
							\end{scope}
						\end{scope}
				\end{scope}
			\end{scope}
			
		\end{scope}
	\end{scope}
}

\pgfplotsset{compat=1.14}

\title{\LARGE \bf Real-Time Quad-Rotor Path Planning Using Convex \\ Optimization and Compound State-Triggered Constraints}

\author{Michael Szmuk$^*$, Danylo Malyuta$^*$, Taylor P. Reynolds$^\dagger$, Margaret Skye Mceowen$^*$, and Beh\c{c}et A\c{c}ikme\c{s}e$^*$
	\thanks{$^*$Autonomous Controls Laboratory, Department of Aeronautics and Astronautics, University of Washington, Seattle, WA 98105, USA {\tt\small \{mszmuk, danylo, skye95, behcet\}@uw.edu}}
	\thanks{$^\dagger$Robotics, Aerospace and Information Networks Laboratory, Department of Aeronautics and Astronautics, University of Washington, Seattle, WA 98105, USA {\tt\small tpr6@uw.edu}}
}

\begin{document}

\maketitle
\thispagestyle{empty}
\pagestyle{empty}


\begin{abstract}

The contribution of this paper is the application of compound state-triggered constraints (STCs) to real-time quad-rotor path planning. Originally developed for rocket landing applications, STCs are made up of a \textit{trigger condition} and a \textit{constraint condition} that are arranged such that satisfaction of the former implies satisfaction of the latter. Compound STCs go a step further by allowing multiple trigger and constraint conditions to be combined via Boolean ``and'' or ``or'' operations. The logical implications embodied by STCs can be formulated using continuous variables, and thus enable the incorporation of discrete decision making into a continuous optimization framework.
In this paper, compound STCs are used to solve quad-rotor path planning problems that would typically require the use of computationally expensive mixed-integer programming techniques. Two scenarios are considered: (1) a quad-rotor flying through a hoop, and (2) a pair of quad-rotors carrying a beam-like payload through an obstacle course. Successive convexification is used to solve the resulting non-convex optimization problem. Monte-Carlo simulation results show that our approach can reliably generate trajectories at rates upwards of 3 and 1.5~Hz for the first and second scenarios, respectively.

\end{abstract}


\section{Introduction} \label{sec:introduction}

The main contribution of this paper is the application of compound state-triggered constraints (STCs) to quad-rotor path planning applications. STCs, and their generalized counterparts, compound STCs, were recently introduced to solve powered-descent guidance rocket landing problems that contained discrete decisions~\cite{szmuk2018jgcdarxiv,szmuk2019scitech,Reynolds2019}. Simply stated, this class of constraints enables user-defined constraint conditions to be enforced if other user-defined trigger conditions are satisfied. To the best of our knowledge, STCs are novel since they capture this discrete logical implication without relying on discrete decision variables. Instead, STCs are formulated using continuous variables, and in practice work well within existing continuous optimization frameworks (e.g. successive convexification). As a result, STCs can be interpreted as \textit{if}-statements that are embedded inside a continuous optimization problem.

Over the past two decades, direct methods for solving optimal control problems have seen a rise in popularity due to the ease of use, performance, and convergence properties offered by modern optimization algorithms~\cite{Betts1998,LiuSurvey2017}. Direct methods are typically used to solve problems with continuous variables, and cannot readily enforce constraints involving discrete decisions. The most common technique used to address this shortcoming is through the use of mixed-integer programming techniques. Despite the existence of efficient branch-and-bound methods, mixed-integer programming techniques suffer from poor computational complexity~\cite{Biegler2014,Richards2015}. The real-time capabilities of such techniques are further hampered when evaluating each set of discrete decisions is expensive.

Mixed-integer programming problems appear in quad-rotor applications quite frequently. In \cite{mip1icra2012}, a centralized mixed-integer quadratic programming (MIQP) approach was used to perform collision avoidance among a team of four heterogeneous quad-rotors. The results showed that feasible solutions could be found in tenths of a second, but that optimality required significantly more computational effort. As discussed in the paper, the methodology was not easily scalable to larger teams of quad-rotors. In \cite{mip2icra2015}, the authors formulated a different MIQP problem to handle the hybrid dynamics of a quad-rotor flying with a mass suspended by a non-rigid string. The paper illustrated that the hybrid nature of the dynamics could be exploited to allow the vehicle to perform otherwise infeasible maneuvers. However, the paper reported computation times upwards of 100~s. In \cite{mip3icra2016}, a Mixed-Integer Semi-Definite Programming (MISDP) approach was proposed to perform aggressive obstacle avoidance in highly cluttered environments (5-26 obstacles). The approach was able to impressively avoid very small obstacles, but reported average computation times of approximately 10 minutes.

In this paper, we propose an STC-based approach that prioritizes computational speed while settling for locally optimal solutions. Two scenarios are used to demonstrate the proposed methodology: (1) a quad-rotor flying through a hoop, and (2) a pair of quad-rotors carrying a beam-like payload through an obstacle course.
The combinatorial elements of these scenarios are formulated into a continuous framework using compound STCs. The successive convexification framework~\cite{mao2016cdc,mao2017aut,szmuk2018scitech} is used to cast the original non-convex problem into a sequence of convex Second-Order Cone Programs (SOCPs).
Our results show that the optimal control problems associated with the first and second scenarios can be solved reliably at average update rates of 15 and 4 Hz, and no slower than 3 and 1.5 Hz, respectively.

In this paper, we adopt the following notation and conventions: an Up-East-North reference frame is used throughout the paper; $\real$, $\real_{+}$, and~$\real_{++}$ are used to denote the set of reals, non-negative reals, and positive reals; $\real^n$, $\real^{m\times n}$, and~$\psd{n}$ are used to denote the space of~$n$-dimensional vectors, $m\times n$-dimensional matrices, and~$n\times n$-dimensional symmetric positive semi-definite matrices; $\set{S}^n\subset\real^{n+1}$ is the unit~$n$-sphere; for vectors quantities, the symbol $\;\hat{}\;$ is used to signify unity norm; $\hat{e}_j$~is used to denote a unit vector with a unity~$\ith{j}{th}$ element; $z\in\real^{\Nz}$~is used to denote a generic solution variable of an optimization problem.

This paper is organized as follows: in~\sref{sec:stcs}, we give a brief overview of STCs and compound STCs; in~\sref{sec:formulation}, we detail our modeling assumptions and the two motivating scenarios; in~\sref{sec:scvx}, we outline the successive convexification algorithm used in the subsequent section; in~\sref{sec:results}, we present our Monte Carlo simulation results for both scenarios; and in~\sref{sec:conclusion}, we provide concluding remarks.


\section{State-Triggered Constraints} \label{sec:stcs}

In this section we provide a concise introduction to STCs. We refer the reader to~\cite{szmuk2018jgcdarxiv,szmuk2019scitech,Reynolds2019} for more details.

\subsection{Logical Statement}

An STC is composed of two parts: a \textit{trigger condition} given by the strict inequality~$g(z)<0$, and a \textit{constraint condition} given by the inequality~$c(z)\leq 0$. We call~$g(z):\real^\Nz\rightarrow\real$ the \textit{trigger function}, and~$c(z):\real^\Nz\rightarrow\real$ the \textit{constraint function}. Both~$g(\cdot)$ and~$c(\cdot)$ are assumed to be differentiable. Formally, an STC enforces the following logical relationship:
\begin{equation} \label{eq:stc_logic}
    g(z) < 0 \;\Rightarrow\; c(z) \leq 0.
\end{equation}
The practical value of an STC is most evident from the contrapositive of~\eqref{eq:stc_logic}, namely, that the constraint condition is not satisfied \textit{only if} the trigger condition is not satisfied.

\subsection{Continuous Formulation}

Mixed-integer programming is the most common framework used to implement discrete decisions such as~\eqref{eq:stc_logic}. However, this approach suffers from poor computational complexity due to the combinatorial nature of integer variables~\cite{Biegler2014,Richards2015}. Moreover, even in the absence of STC-like constraints, practical (non-convex) path planning problems often require the use of sequential solution methods (e.g. Sequential Quadratic Programming, Successive Convexification). For these reasons we seek a continuous formulation of~\eqref{eq:stc_logic} that is amenable to a sequential (continuous) implementation without incurring the added computational complexity of mixed-integer approaches.

An equivalent continuous formulation of~\eqref{eq:stc_logic} was introduced in~\cite{szmuk2018jgcdarxiv}, and is given by
\begin{equation} \label{eq:stc_cont}
    h(z) \definedas \shat(z) \cdot c(z) \leq 0,
\end{equation}
where~$\shat(z)\definedas -\min\big(0,g(z)\big)$. By inspection, we see that if~$g(z)<0$, then~$\shat(z)>0$, and~\eqref{eq:stc_cont} reduces to~$c(z)\leq 0$. In contrast, if~$g(z)\geq 0$, then~$\shat(z)=0$, and~\eqref{eq:stc_cont} is trivially satisfied for any value of~$c(z)$ (i.e. the constraint condition is not enforced). Thus, we conclude that~\eqref{eq:stc_logic} and~\eqref{eq:stc_cont} are logically equivalent, and emphasize that the latter can be implemented in a continuous optimization framework.

\subsection{Compound State-Triggered Constraints}

Compound STCs were introduced in~\cite{szmuk2019scitech}, and are a generalization of the scalar STC formulation given in~\eqref{eq:stc_logic} and~\eqref{eq:stc_cont}. Compound STCs have trigger and constraint conditions that are composed using Boolean ``and'' or ``or'' operations. Here, we present compound STCs with ``and''- and ``or''-trigger conditions, and ``or''-constraint conditions. The logical representations of these STCs are given by
\begin{subequations} \label{eq:cstc_logic}
    \begin{align}
        \bigwedge_{j=1}^{\Ng}\big(g_j(z) < 0\big)\;&\Rightarrow\;\bigvee_{j=1}^{\Nc}\big(c_j(z)\leq 0\big),\label{eq:cstc_and_logic} \\
        \bigvee_{  j=1}^{\Ng}\big(g_j(z) < 0\big)\;&\Rightarrow\;\bigvee_{j=1}^{\Nc}\big(c_j(z)\leq 0\big),\label{eq:cstc_or_logic}
    \end{align}
\end{subequations}
where there are~$\Ng$ trigger conditions, $\Nc$ constraint conditions, and each~$g_j(\cdot)$ and~$c_j(\cdot)$ is defined as in the scalar case. The corresponding continuous formulations are given by
\begin{subequations} \label{eq:cstc_cont}
    \begin{align}
        h_{\land}(z) &\definedas \Bigg[\prod_{j=1}^{\Ng}\hspace{0.05cm}\shat_j(z)\Bigg]\cdot \Bigg[\prod_{j=1}^{\Nc} \big(c_j(z)+\alpha_j\big) \Bigg] = 0, \label{eq:cstc_and_cont} \\
        h_{\lor }(z) &\definedas \Bigg[\sum_{ j=1}^{\Ng}\shat_j(z)\Bigg]\cdot \Bigg[\prod_{j=1}^{\Nc} \big(c_j(z)+\alpha_j\big) \Bigg] = 0, \label{eq:cstc_or_cont}
    \end{align}
\end{subequations}
where~$\alpha_j\in\real_{+}$ are non-negative slack variables, and each $\shat_j(\cdot)$ is defined as in the scalar case.

We conclude this section with two comments. First, formulations with equality constraint conditions can be obtained by substituting equalities in place of the (non-strict) inequalities in~\eqref{eq:stc_logic}-\eqref{eq:cstc_logic}, and omitting the slack variables in~\eqref{eq:cstc_cont}. Second, a compound STC with an ``and''-constraint condition is emulated by enforcing~$\Nc$ separate compound STCs, each with the original compound trigger condition and one of the (scalar) constraint conditions. We are now ready to apply STCs to the scenarios detailed in the next section.


\section{Problem Formulation} \label{sec:formulation}

\begin{figure*}[t]
    \centering
    \input{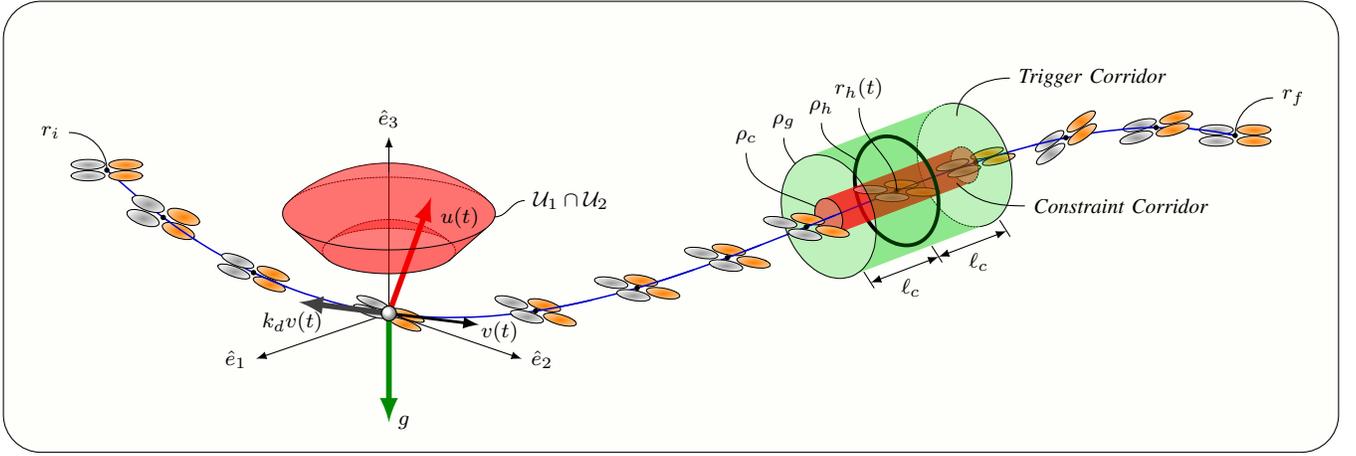}
    \vspace{-0.5cm}
    \caption{Illustration of Scenario 1: The feasible thrust set is shown as the red volume on the left. The velocity, thrust, drag, and gravity vectors are shown as the black, red, gray, and green vectors, respectively. The initial and final positions are indicated on the left and right, respectively. The hoop, trigger corridor, and constraint corridor are shown on the right. The vector~$\nh$ is orthogonal to the plane of the hoop, and is not shown. The compound STC makes the volume that is inside the trigger corridor but outside of the constraint corridor infeasible.} 
    \vspace{-0.5cm}
    \label{fig:scenario1}
\end{figure*}

In this section, we outline two quad-rotor path planning scenarios using compound STCs. Before presenting the two scenarios, we briefly discuss our assumed guidance and control (G\&C) architecture and simplified quad-rotor model.

\subsection{G\&C Architecture}

The G\&C architecture assumed in this paper separates guidance and control into two distinct tasks. The \textit{guidance task} involves generating an open-loop trajectory at a low frequency, whereas the \textit{control task} involves computing high frequency closed-loop control actions to stay on the guidance trajectory. The objective of the guidance task is to ensure feasibility (e.g. respecting vehicle dynamics and control limits, and avoiding obstacles), while the objective of the control task is to provide robustness to plant uncertainties and external disturbances (e.g. wind gusts, and battery voltage variability). The control task is typically subdivided into hierarchically arranged thrust, attitude, and translation controllers.

\subsection{Simplified Quad-Rotor Dynamics} \label{sec:dynamics}

We assume a 3-DoF quad-rotor dynamics model, similar to the one used in~\cite{szmuk2017iros,szmuk2018iros}. This model is given by
\begin{gather*}
    \dot{x}(t) = A x(t) + B u(t) + E w, \\[1.0ex]
    A \definedas \mat{0_{3\times 3} & \hspace{0.61cm}I_{3\times 3} \\ 0_{3\times 3} & -\kd I_{3\times 3}},\; B \definedas \frac{1}{m}\mat{0_{3\times 3} \\ I_{3\times 3}}, \\[1.0ex]
    x(t)\definedas\mat{r^\transp(t) & v^\transp(t)}^\transp,\; E \definedas -\hat{e}_4,\; w \definedas g,
\end{gather*}
where~$r(t)\in\real^3$ is the position state, $v(t)\in\real^3$ is the velocity state, $u(t)\in\set{U}\subset\real^3$ is the thrust (control) vector, $m\in\real_{++}$~is the mass of the vehicle, $\kd\in\real_{+}$~is the drag coefficient, and $g\in\real_{++}$~is the local gravitational acceleration. Note that the drag model is simplified due to its linear dependence on~$v(t)$.

To define the control set~$\set{U}$, we first define three sets. The first set represents the allowable thrust magnitudes of the vehicle, and is given by
\begin{equation*}
    \set{U}_1\definedas\big\{u\in\real^3 : 0<\Tmin\leq\|u\|_2\leq\Tmax\big\},
\end{equation*}
where~$\Tmin$ and~$\Tmax$ are the minimum and maximum allowable thrust magnitudes. In practice, these bounds are selected conservatively to ensure that the underlying controllers can command thrust and torques independently. Note that~$\set{U}_1$ is non-convex. The second set represents the allowable tilt angles of the vehicle, and is given by
\begin{equation*}
    \set{U}_2\definedas\big\{u\in\real^3 : \cos\thetamax\|u\|_2\leq\hat{e}_1^\transp u\big\},
\end{equation*}
where~$\thetamax\in(0\degree,180\degree)$ is the maximum allowable tilt angle. Note that~$\set{U}_2$ is non-convex for~$\thetamax>90\degree$. The third set represents the thrust vectors with a vertical component equal and opposite to the weight of the vehicle (i.e. control inputs that maintain a constant altitude). This set is given by
\begin{equation*}
    \set{U}_3\definedas\big\{u\in\real^3 : \hat{e}_1^\transp u=mg\big\},
\end{equation*}
and has a non-empty interior when~$\Tmin \leq mg < \Tmax$.

For three-dimensional applications,~$\set{U}=\set{U}_1 \cap \set{U}_2$ is non-convex, and the optimal control problem can be convexified using the lossless convexification technique introduced in~\cite{behcet2007jgcd} (also see~\cite{szmuk2017iros}). For two-dimensional applications requiring only horizontal motion,~$\set{U}=\set{U}_1 \cap \set{U}_2 \cap \set{U}_3$ is convex, and the vertical dimension of the problem can be omitted from the formulation of the guidance problem. 

\subsection{Scenario 1: Quad-Rotor Flying Through a Hoop} \label{sec:scen1}

\begin{figure*}[t!]
    \centering
    \input{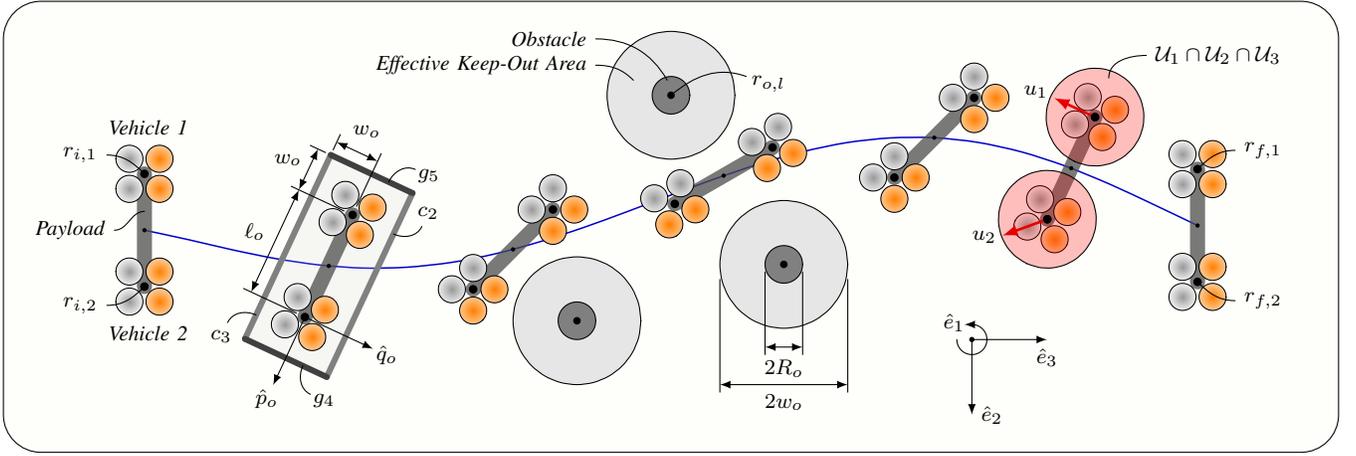}
    \vspace{-0.5cm}
    \caption{Illustration of Scenario 2: The two vehicles and the payload are shown in their initial state on the left. The keep-out region defined by~\eqref{eq:ex2_cstc_logic} and~\eqref{eq:ex2_cstc_cont} is shown at the second time instance. This region effectively keeps the vehicles and the payload outside of the \textit{effective keep-out areas} shown around the obstacles. Lastly, the convex control sets are represented by the red circles in the second to last time instance.}
    \vspace{-0.5cm}
    \label{fig:scenario2}
\end{figure*}

The first scenario involves flying a quad-rotor through a hoop, and is three dimensional in nature (see Figure~\ref{fig:scenario1}). We formulate this scenario as a fixed-final-time optimal control problem of duration~$\tf\in\real_{++}$.

The boundary conditions of this problem are prescribed through the initial and final position vectors~$r_i$ and~$r_f$. The initial and final velocities are assumed to equal zero, and the initial and final controls are assumed to be equal and opposite the weight of the vehicle.

The position of the center of the hoop is denoted by $\rh(t)\in\real^3$, and can vary with time. We assume that~$\rh(t)$ is known for all~$t\in[0,\,\tf]$. We represent the hoop's orientation using an orthogonal unit vector~$\nh\in\set{S}^2$, which we assume is constant for simplicity. We note that our formulation readily handles small variations in hoop orientation, and emphasize that~$\rh(t)$ and~$\nh$ \textit{are not} decision variables of the optimal control problem.

\boxing{b!}{0.457}{-0.3cm}{-0.25cm}{
	\begin{problem}{\bf \textit{Non-Convex Formulation of Scenario 1}} \label{problem:scenario1}
		$$ \underset{u}{\text{minimize}}\;\int_{0}^{t_f}{\|u(t)\|_2 dt}$$
		\hspace{0.25cm} subject to: \vspace{-0.15cm}
		\begin{gather*}
			r(0) = r_i,\, r(t_f) = r_f,\\ v(0)=v(t_f) = 0_{3\times 1},\\ u(0) = u(t_f) = mg\hat{e}_1,\\
			\dot{x}(t) = Ax(t)+Bu(t)+Ew,\, u(t)\in\set{U}_1\cap\set{U}_2,\\ \|v(t)\|_2 \leq \vmax,\, h_1\big(r(t),\rh(t)\big) \leq 0.
		\end{gather*}
	\end{problem}
}

The constraint imposed by the hoop is formulated using a compound STC that restricts the path of the quad-rotor to a \textit{constraint corridor} whenever the vehicle is inside a \textit{trigger corridor}. The constraint corridor functions to guide the vehicle through the hoop without collision. The geometry of this constraint is illustrated in Figure~\ref{fig:scenario1}, where~$\Rh$ denotes the radius of the hoop. Omitting the time arguments of~$r(t)$ and~$\rh(t)$, the logical implication of this compound STC is given by

\begin{subequations} \label{eq:ex1_cstc_logic}
    \begin{equation}
        \bigwedge_{j=1}^{3} \big(g_j(r,\rh) < 0\big) \;\Rightarrow\; \big(c_1(r,\rh) \leq 0\big),
    \end{equation}
    \begin{align}
        g_1(r,\rh) &\definedas \nh^\transp(\rh-r)-\Lh, \\
        g_2(r,\rh) &\definedas \nh^\transp(r-\rh)-\Lh, \\
        g_3(r,\rh) &\definedas (r-\rh)^\transp\Nho^\transp\Nho(r-\rh)-\Rhg^2, \\
        c_1(r,\rh) &\definedas (r-\rh)^\transp\Nho^\transp\Nho(r-\rh)-\Rhc^2,
    \end{align}
\end{subequations}
where $\Lh\in\real_{++}$ is the half-length of the corridor, $\Rhg$ is the radius of the trigger corridor, $\Rhc$ is the radius of the constraint corridor, and~$\Nho\definedas I_{3\times 3}-\nh\nh^\transp$. In practice, these parameters are selected to satisfy~$0 \leq \Rhc \ll \Rh \ll \Rhg$. From~\eqref{eq:cstc_and_cont} and~\eqref{eq:ex1_cstc_logic}, we obtain the following continuous formulation:
\begin{equation} \label{eq:ex1_cstc_form}
    h_1(r,\rh) \definedas \Bigg[\prod_{j=1}^{3}\shat_j(r,\rh)\Bigg]\cdot c_1(r,\rh) \leq 0,
\end{equation}
where~$\shat_j(r,\rh)\definedas -\min\big(0,g_j(r,\rh)\big)$. Since~$\Nc=1$ in this case, we omit the slack variable from~\eqref{eq:cstc_and_cont}, and replace the equality with an inequality. The associated fuel-optimal non-convex optimal control problem is summarized in Problem~\ref{problem:scenario1}. This problem has two sources of non-convexity: the control set, and the compound STC given in~\eqref{eq:ex1_cstc_form}. The former is addressed using lossless convexification (see~\sref{sec:dynamics}), whereas the latter is convexified using successive convexification (see~\sref{sec:scvx}).

We conclude this section with a few remarks. First, note that Problem~\ref{problem:scenario1} includes a constraint that limits the velocity to a maximum of~$\vmax\in\real_{++}$. This constraint is added in order to mitigate constraint clipping introduced by temporal discretization (see~\sref{sec:results}). Second, the above scenario is similar to the \textit{Agile Flip Maneuver} presented in~\cite{szmuk2017iros}, which required the quad-rotor to maneuver through a waypoint defined midway along the trajectory. However, the key difference between~\cite{szmuk2017iros} and the scenario described in this section is that the compound STC enables the optimization to choose \textit{if} and \textit{when} the trajectory will pass through the hoop. Third, the above formulation can be modified such that the vehicle is \textit{required} to pass through the hoop. This can be done either by selecting~$\Rhg$ sufficiently large, or by omitting~$g_3(\cdot)$ from the formulation. Lastly, the direction and speed of the trajectory at the hoop can be specified by enforcing an additional compound STC with an appropriate velocity-dependent constraint condition.


\subsection{Scenario 2: Cooperative Obstacle Avoidance} \label{sec:scen2}

The second scenario consists of two identical quad-rotors cooperatively negotiating an obstacle course (see Figure~\ref{fig:scenario2}). We restrict the motion of the vehicles to the horizontal plane, hence making this a two-dimensional scenario. As in the first scenario, we treat this problem as a fixed-final-time problem, and use the subscripts~1 and~2 to distinguish between quantities associated with the two vehicles.

The position boundary conditions are given by the initial position vectors~$r_{i,1}$ and~$r_{i,2}$, and the final position vectors~$r_{f,1}$ and~$r_{f,2}$. For each vehicle, the velocity and control boundary conditions are identical to those used in the first scenario.

The quad-rotors are linked together by a beam-like payload of length~$\Lspacing\in\real_{++}$, modeled by the following non-convex equality constraint:
\begin{equation} \label{eq:ex2_fixed_dist}
    \|r_1(t)-r_2(t)\|_2 = \Lspacing.
\end{equation}
We assume that the vehicles maintain their ability to control their attitudes independently of one another (i.e. each vehicle can control its attitude as in~\sref{sec:scen1}), and that the boundary conditions are feasible with respect to~\eqref{eq:ex2_fixed_dist}.

The flight space contains~$\No$ stationary cylindrical obstacles of identical radius~$\Ro\in\real_{++}$. The position of each obstacle~$\jo\in\obset\definedas\{1,\ldots,\No\}$ is denoted by~$\roj\in\real^3$.
Each obstacle is assumed to vertically span the available space.

This scenario is challenging since the absence of vehicle-obstacle collisions does not guarantee the absence of payload-obstacle collisions. To address this issue, a compound STC is used to define a keep-out rectangle around the two vehicles. Omitting the time arguments of~$r_1(t)$ and~$r_2(t)$, the logical implication of this compound STC is given by
\begin{subequations} \label{eq:ex2_cstc_logic}
    \begin{equation}
        \bigwedge_{j=1}^{2}\big(g_{j+3}(\cdot,\cdot,\cdot)<0\big)\;\Rightarrow\;\bigvee_{j=1}^{2}\big(c_{j+1}(\cdot,\cdot,\cdot)\leq 0\big),
    \end{equation}\vspace{-0.40cm}
    \begin{align}
        g_4(r_1,r_2,\roj) &\definedas \phato^\transp(r_1-\roj)-\Wo,\\
        g_5(r_1,r_2,\roj) &\definedas \phato^\transp(\roj-r_2)-\Wo,\\
        c_2(r_1,r_2,\roj) &\definedas \qhato^\transp(\roj-r_2)+\Wo,\\
        c_3(r_1,r_2,\roj) &\definedas \qhato^\transp(r_1-\roj)+\Wo,
    \end{align}
\end{subequations}
where~$\phato\definedas(r_2-r_1)/\|r_2-r_1\|_2$, $\qhato$ is orthogonal to~$\phato$, and~$\Wo\in\real_{++}$ is the minimum spacing enforced around each vehicle. In practice, $\Wo$ is selected such that~$\Wo\geq\Ro$. These quantities are illustrated in Figure~\ref{fig:scenario2}. From~\eqref{eq:cstc_and_cont} and~\eqref{eq:ex2_cstc_logic}, we obtain the following continuous formulation:
\begin{equation} \label{eq:ex2_cstc_cont}
    \begin{split}
    h_2(r_1,r_2,\roj)\definedas
         &\Bigg[\prod_{j=1}^{2}\shat_{j+3}(r_1,r_2,\roj)\Bigg]\\
    \cdot&\Bigg[\prod_{j=1}^{2}\big(c_{j+1}(r_1,r_2,\roj)+\alpha_j\big)\Bigg]=0,
    \end{split}
\end{equation}
where~$\shat_j(\cdot,\cdot,\cdot)\definedas -\min\big(0,g_j(\cdot,\cdot,\cdot)\big)$, and $\alpha_1,\alpha_2\in\real_{+}$ are non-negative slack variables as in~\eqref{eq:cstc_cont}. Defining the following quantities
\begin{gather*}
    \tilde{A}\definedas\blkdiag{A, A},\, \tilde{B}\definedas\blkdiag{B, B},\\ \tilde{E}\definedas\blkdiag{E, E},\, \tilde{w}\definedas[w \;\, w]^\transp,\\ \tilde{x}\definedas [x_1^\transp \;\, x_2^\transp]^\transp,\,
    \tilde{u}\definedas [u_1^\transp \;\, u_2^\transp]^\transp,
\end{gather*}
the associated fuel-optimal non-convex optimal control problem is summarized in Problem~\ref{problem:scenario2}. This problem has two sources of non-convexity: the equality constraint given in~\eqref{eq:ex2_fixed_dist}, and the compound STC given in~\eqref{eq:ex2_cstc_cont}. Both of these non-convexities are handled using successive convexification (see~\sref{sec:scvx}).

\boxing{t!}{0.457}{0.20cm}{-0.7cm}{
	\begin{problem}{\bf \textit{Non-Convex Formulation of Scenario 2}} \label{problem:scenario2}
		$$ \underset{u}{\text{minimize}}\;\int_{0}^{t_f}{\big(\|u_1(t)\|_2+\|u_2(t)\|_2\big) dt}$$
		\hspace{0.25cm} subject to: \vspace{-0.15cm}
		    \begin{gather*}
		        r_1(0) = r_{i,1},\, r_1(t_f) = r_{f,1},\\
		        r_2(0) = r_{i,2},\, r_2(t_f) = r_{f,2},\\
		        v_1(0) = v_2(0) = v_1(t_f) = v_2(t_f) = 0_{3\times 1},\\
		        u_1(0) = u_2(0) = u_1(t_f) = u_2(t_f) = mg\hat{e}_1,\\[1.5ex]
			    \dot{\tilde{x}}(t) = \tilde{A}\tilde{x}(t)+\tilde{B}\tilde{u}(t)+\tilde{E}\tilde{w},\\
			    u_1(t),u_2(t)\in\set{U}_1\cap\set{U}_2\cap\set{U}_3, \\[1.5ex]
			    h_2\big(r_1(t),r_2(t),\roj\big)=0,\;\;\forall j\in\obset, \\
			    \|v_1(t)\|_2\leq \vmax,\,\|v_2(t)\|_2\leq\vmax,\\ \|r_1(t)-r_2(t)\|_2 = \Lspacing.
		    \end{gather*}
	\end{problem}
}

We conclude this section with three remarks. First, the control set is convex due to the intersection of~$\set{U}_1\cap\set{U}_2$ with~$\set{U}_3$. Second, implementations of Problem~\ref{problem:scenario2} can omit the vertical dimension of the problem due to the two-dimensional nature of this scenario. Third, the compound STC in~\eqref{eq:ex2_cstc_logic} and~\eqref{eq:ex2_cstc_cont} can be extended to more complicated geometries (e.g. multiple vehicles vehicles carrying an \textit{L}-shaped payload).


\section{Successive Convexification} \label{sec:scvx}

In this section, we provide an overview of the successive convexification algorithm used to solve Problems~\ref{problem:scenario1} and~\ref{problem:scenario2}. We refer the reader to~\cite{szmuk2018jgcdarxiv} for more details.

Successive convexification is a framework that solves non-convex continuous-time optimal control problems by solving a sequence of convex discrete-time parameter optimization \textit{subproblems}. Each subproblem is an SOCP that approximates the original problem by linearizing non-convexities about the previous iteration, and is obtained using two steps: discretization and linearization.

\subsection{Discretization \& Linearization} \label{sec:disc_lin}

The (temporal) \textit{discretization step} divides the time horizon of the optimal control problem into~$\KK-1$ temporal intervals of length~$\dt\definedas\tf/(\KK-1)$. For each node~$k\in\Kset\definedas\{1,2,\ldots,\KK\}$, the time is given by~$\tk\definedas (k-1)\dt$.

Since the dynamics of Problems~\ref{problem:scenario1} and~\ref{problem:scenario2} are linear time-invariant, the discrete-time dynamics can be expressed analytically as a function of~$\tf$. Assuming a first-order-hold on the control, we represent these discrete-time dynamics by the quantities~$\xdk{k}\in\real^{\Nx}$, $\udk{k}\in\real^{\Nu}$, $\Ad\in\real^{\Nx\times\Nx}$, $\Bd^{-},\Bd^{+}\in\real^{\Nx\times\Nu}$, $\Ed\in\real^{\Nx\times\Nw}$, and~$\Wd\in\real^{\Nw}$. The discrete-time state and control constraints are obtained by enforcing said constraints at each temporal node.

The \textit{linearization step} linearizes the non-convexities that cannot be convexified using lossless convexification (i.e.~\eqref{eq:ex1_cstc_form}, \eqref{eq:ex2_fixed_dist}, and~\eqref{eq:ex2_cstc_cont}). Since this approximation is only made to first order, the subproblem is \textit{guaranteed} to be convex. However, the linearization also introduces two issues: \textit{artificial infeasibility} and \textit{artificial unboundedness}. 

To aide in the ensuing explanation, we define~$\Ksetm\definedas\Kset\setminus\KK$, $\ud\definedas[\udk{1}^\transp,\ldots,\udk{\KK}^\transp]^\transp$, $\zdk{k}\definedas[\xdk{k}^\transp,\,\udk{k}^\transp,\,\adk{k}]^\transp$, and~$\zd\definedas[\zdk{1}^\transp,\ldots,\zdk{\KK}^\transp]^\transp$, where~$\adk{k}\in\real_{+}^{\Na}$ is a vector of non-negative slack variables and~$\Nz=\Nx+\Nu+\Na$. We concatenate the non-convex state constraints into the vector-valued equality constraint~$h(\zdk{k}) = 0$.

Artificial infeasibility is resolved by adding \textit{virtual control} terms~$\vd\in\real^{\Nx}$ to the discrete-time dynamics, and augmenting the cost with~$\Jvc(\vd)\definedas\sum_{k\in\Ksetm}\|\Wvc\vdk{k}\|_1$, where~$\Wvc\in\psd{\Nx}$ is a user-specified weight matrix, and $\vd\definedas[\vdk{1}^\transp,\ldots,\vdk{\KK-1}^\transp]^\transp$.
The addition of~$\Jvc(\cdot)$ penalizes violations of the dynamics, and allows dynamic infeasibility to occur (if necessary) during the convergence process.

Artificial unboundedness is resolved by augmenting the cost with~$\Jtr(\zd)\definedas \sum_{k\in\Kset}\delta\zdk{k}^\transp\Wtr\delta\zdk{k}$, where~$\Wtr\in\psd{\Nz}$ is a user-specified weight matrix, $\delta\zdk{k}\definedas\zdk{k}-\zdk{k}^{*}$, and~$\zdk{k}^{*}$ denotes the solution obtained during the previous iteration. The addition of~$\Jtr(\cdot)$ ensures that Problem~\ref{problem:socp} remains bounded, and keeps the solution close to the linearization point.

\subsection{Subproblem}

\boxing{t!}{0.457}{0.18cm}{-0.7cm}{
	\begin{problem}{\bf \textit{Convex Subproblem (SOCP)}} \label{problem:socp}
		$$ \underset{\ud,\vd}{\text{minimize}}\;\Jorig(\zd)+\Jtr(\zd)+\Jvc(\vd)$$
		\hspace{0.25cm} subject to: \vspace{-0.15cm}
		    \begin{gather*}
			    \xdk{1} = \xdk{d,i},\,\xdk{\KK} = \xdk{d,f},\,\udk{1} = \udk{d,i},\,\udk{\KK} = \udk{d,f},\\
			    \begin{split}
			        \xdk{k+1} = \Ad\xdk{k}+&\Bd^{-}\udk{k}+\Bd^{+}\udk{k+1} \\ +&\Ed\Wd+\vdk{k},\;\;\forall k\in\Ksetm,
			    \end{split} \\[0.0ex]
    			h(\zdk{k}^{*})+\left.\frac{\partial h}{\partial \zdk{k}}\right|_{\zdk{k}^*}\hspace{-0.2cm} \delta\zdk{k}=0,\;\udk{k}\in\set{U},\;\;\forall k\in\Kset.
            \end{gather*}
	\end{problem}
}

Problem~\ref{problem:socp} summarizes the subproblem used in the proposed lossless convexification algorithm. This problem consists of (\hspace{0.01cm}\textit{i}\hspace{0.02cm}) an objective function made up of the original objective~$\Jorig(\zd)$ and the two augmented cost terms discussed in~\sref{sec:disc_lin}; (\hspace{0.01cm}\textit{ii}\hspace{0.02cm}) the boundary conditions denoted by~$\xdk{d,i}$, $\udk{d,i}$, $\xdk{d,f}$, and~$\udk{d,f}$; (\hspace{0.01cm}\textit{iii}\hspace{0.02cm}) the discrete-time dynamics with virtual control; (\hspace{0.01cm}\textit{iv}\hspace{0.02cm}) the non-convex state constraints linearized about the previous solution; and (\hspace{0.01cm}\textit{v}\hspace{0.02cm}) the control set~$\set{U}$ detailed in~\sref{sec:dynamics}. As stated previously, we assume that non-convexity in~$\set{U}$ is handled by lossless convexification~\cite{szmuk2017iros}.

\subsection{Algorithm}

The successive convexification algorithm is outlined in Algorithm~\ref{algorithm:scvx}, and is a simplified version of the soft-trust-region algorithm presented in~\cite{szmuk2018jgcdarxiv}. The algorithm is initialized by user-specified problem data (e.g.~$\tf$, boundary conditions) and an initialization trajectory generated by linearly interpolating between the specified boundary conditions. The subsequent discretization step computes the discrete-time dynamics. The algorithm then enters a loop that successively linearizes and solves Problem~\ref{problem:socp} until the trust region and virtual control costs~$\Jtr(\cdot)$ and~$\Jvc(\cdot)$ are less than their respective specified thresholds $\trtol,\vctol\in\real_{++}$. Upon convergence, the algorithm returns the converged solution~$\zd$. If the algorithm does not converge within a set number of iterations, then~$\tf$ is increased and the problem is resolved.

\algrenewcommand\alglinenumber[1]{\normalsize #1:}
\newcommand{\algemph}[1]{\textbf{\textit{#1}}}
\newcommand{\algvar}[1]{\textit{\texttt{#1}}}
\makeatletter
\xpatchcmd{\algorithmic}{\itemsep\z@}{\itemsep=0.05cm}{}{}
\makeatother
\boxing{t!}{0.457}{0.20cm}{-0.7cm}{
	\begin{algorithm}{\bf \textit{Successive Convexification (Soft-TR)}} \label{algorithm:scvx}
	    \vspace{0.1cm}
		\begin{algorithmic}[1] \normalsize
		    \State \algemph{initialize} - provide problem data, and~$\zdk{k}^{*}$ for all $k\in\Kset$
		    \State \algemph{discretize} - compute~$\Ad$, $\Bd^{-}$, $\Bd^{+}$, $\Ed$, $\Wd$
		    \State \algemph{set} $\algvar{converged}=0$
		    \While{$(\algvar{converged} = 0)$}
		        \State \algemph{linearize} - compute $h(\zdk{k}^{*})$ and~$\left.\partial h/\partial\zdk{k}\right|_{\zdk{k}^{*}}$
    		    \State \algemph{solve SOCP} - compute solution for Problem~\ref{problem:socp}
    		    \If{$(\Jtr(\zd)<\trtol)\land(\Jvc(\vd)<\vctol)$}
    		        \State {$\algvar{converged}=1$}
    		    \Else
    		        \State {$\zdk{k}^{*} \gets \zdk{k}$ for all $k\in\Kset$}
    		    \EndIf
		    \EndWhile
		    \State \algemph{return} - computed solution $\zdk{k}$ for all $k\in\Kset$
		\end{algorithmic}
	\end{algorithm}
}

\section{Results} \label{sec:results}

This section presents Monte-Carlo simulation results for the two scenarios presented in
\sref{sec:formulation}. We focus on measuring the runtime of
Algorithm~\ref{algorithm:scvx} and quantifying two primary failure modes of Algorithm~\ref{algorithm:scvx}:
failure to converge in less than 20 iterations, and inter-sample constraint violation. All
results were obtained on a desktop PC running a Ubuntu 18.04.1 operating system with a 3.60~GHz Intel Cote i7-6850K processor and 64~GB of RAM. MATLAB was used to run the ECOS~\cite{ecos} solver through the CVX~\cite{cvx} parsing interface, and timing data were obtained from the \textit{solve time} parameter returned by the \texttt{cvx\_toc} function.



\subsection{Scenario 1} \label{sec:results_secn1}

\begin{figure}[t!]
    \centering
    \includegraphics[width=0.457\textwidth]{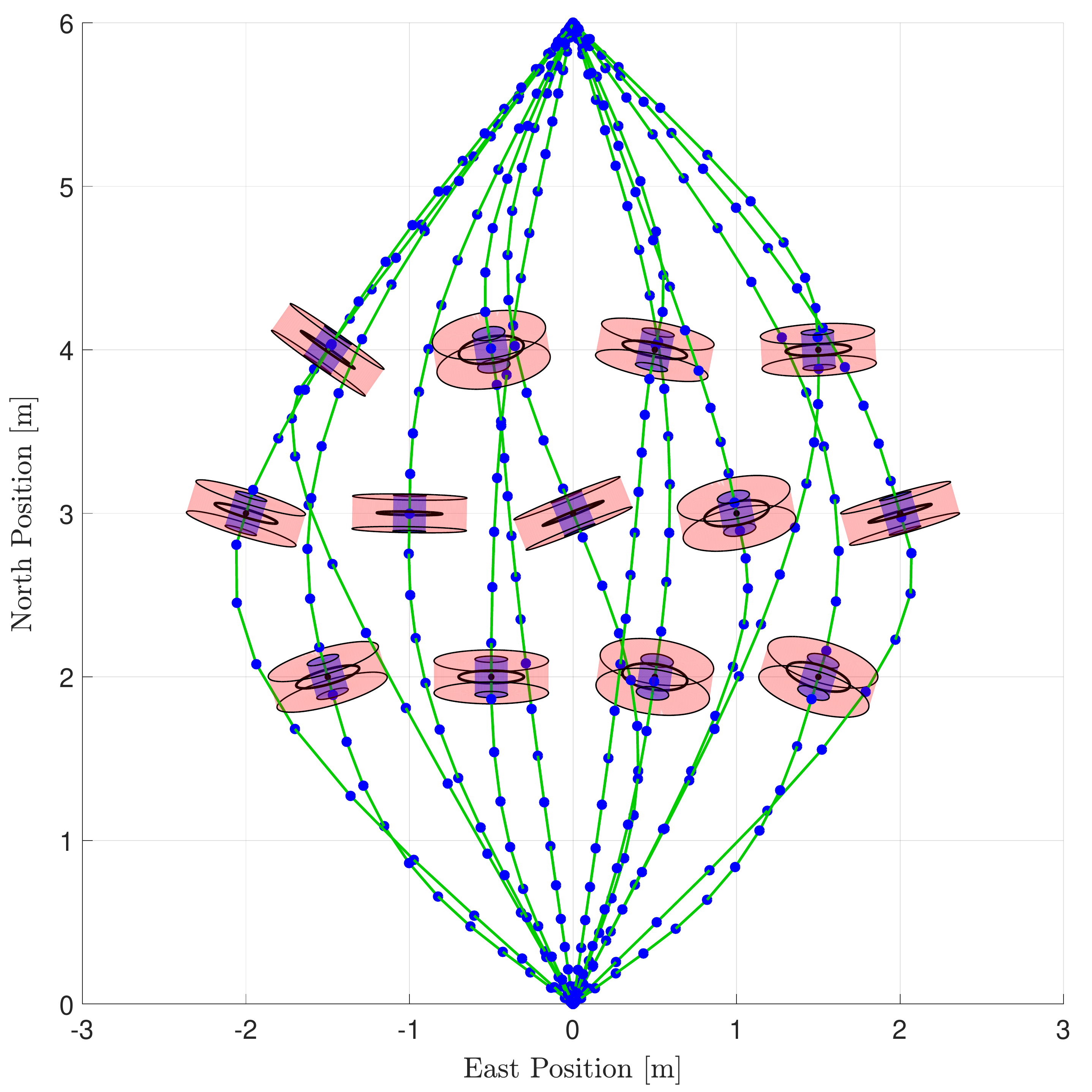}
    \caption{Scenario 1 sample trajectories for $\KK=30$ using a grid of hoop positions and randomized hoop orientations. Vehicle motion is from top to bottom. The blue dots represent temporal nodes, and the trigger and constraint corridors are not drawn to scale. Each trajectory is constrained to pass through one of the hoops by omitting~$g_3(\cdot)$ from~\eqref{eq:ex1_cstc_logic}-\eqref{eq:ex1_cstc_form}.}
    \vspace{-0.6cm}
    \label{fig:scen1_results}
\end{figure}

For this scenario, the performance of Algorithm~\ref{algorithm:scvx} was evaluated by running multiple cases with identical boundary conditions and final time, but with different hoop locations and orientations. Example simulations are presented in Figure~\ref{fig:scen1_results}, which shows 13 trajectories obtained using a temporal resolution of~$\KK=30$, randomized hoop orientations, and a manually selected grid of hoop positions. The hoop positions, tilt angles, and heading angles were sampled from uniform random distributions given by $\rh\in [-1,1]\times [-2,2]\times [2,4]$~m, $\phi\in [-25\degree,25\degree]$, and $\psi\in [-35\degree,35\degree]$, respectively. Other problem parameters were set to the following constant values:
\begin{gather*}
  \tf = 4~\text{s},\quad
  ,\quad r_i=(0,0,0)~\text{m}, \\
  r_f=(0,0,6)~\text{m},\quad \Lh=0.5~\text{m},\quad \Rhc = 0~\text{m}, \\
  \Rhg=\infty~\text{m},\quad
  \vmax=2\Lh/\dt~\text{m/s},\quad \Tmin=2~\text{N}, \\
  \Tmax=5~\text{N},\quad\thetamax=45\degree,\quad m=0.35~\text{kg}, \\
  g=9.81~\text{m/s}^2,\quad \kd=0~\text{s}^{-1},\quad
  \Wvc = 10^5\cdot I_{6\times 6}, \\
  \Wtr = 0.1\cdot\blkdiag{I_{3\times 3},0_{6\times 6}}.
\end{gather*}

Monte Carlo batch runs were conducted for temporal resolutions $\KK=15,20,25,30$. Each batch consisted of 100 cases. Cases that resulted in infeasibility were noted, and new cases were sampled in their place. Tables~\ref{table:scen1_timing} and \ref{table:scen1_clipping} provide statistics for Algorithm~\ref{algorithm:scvx} runtime and inter-sample
constraint violation. As expected, the runtime increased with $\KK$ since the number of decision variables in Problem~\ref{problem:socp} increased. Nevertheless, the algorithm demonstrated the ability to run at interactive
rates even for $\KK=30$, with a worst-case observed execution rate of approximately $3~\text{Hz}$.

\begin{table}[b!]
    \begin{center}
        \vspace{-0.3cm}
    	\caption{Algorithm~\ref{algorithm:scvx} runtime for Scenario 1 [ms].}
 		\begin{tabu} to 0.457\textwidth {c|X[c]X[c]X[c]X[c]X[c]}
 			\hhline{======} \\[-1.50ex]
     		\KK  & Mean & Median & Std. Dev. & Min & Max \\[0.5ex]
     		\hhline{------} \\[-1.50ex]
     		$15$ & $34$ & $31$ & $25$ & $15$ & $151$ \\[1ex]
            $20$ & $43$ & $34$ & $32$ & $19$ & $219$ \\[1ex]
            $25$ & $55$ & $42$ & $47$ & $24$ & $306$ \\[1ex]
            $30$ & $67$ & $52$ & $53$ & $29$ & $337$ \\[1ex]
     		\hhline{======}
 		\end{tabu}
        \label{table:scen1_timing}
        \vspace{-0.3cm}
 	\end{center}
\end{table}
\begin{table}[b!]
    \begin{center}
    	\caption{Constraint clipping for Scenario 1 [cm].}
 		\begin{tabu} to 0.457\textwidth {c|X[c]X[c]X[c]X[c]X[c]}
 			\hhline{======} \\[-1.50ex]
     		\KK  & Mean & Median & Std. Dev. & Min & Max \\[0.5ex]
     		\hhline{------} \\[-1.50ex]
     		$15$ & $37$ & $32$ & $37$ & $2$ & $224$ \\[1ex]
            $20$ & $15$ & $10$ & $14$ & $1$ & $100$ \\[1ex]
            $25$ & $6$ & $6$ & $3$ & $1$ & $16$ \\[1ex]
            $30$ & $4$ & $4$ & $2$ & $0$ & $10$ \\[1ex]
     		\hhline{======}
 		\end{tabu}
        \label{table:scen1_clipping}
        \vspace{-0.2cm}
 	\end{center}
\end{table}

Additionally, it was observed that inter-sample constraint violation decreased as $\KK$
increased. This was expected, since a smaller $\dt$ reduces the possibility
of the quad-rotor ``hopping'' over the hoop in one discrete time step.

Lastly, our simulations resulted in 400 converged cases, and 19 failures to converge. In each case, the cause of infeasibility was an excessively short final time, and feasibility was recovered by selecting a larger~$\tf$.

\subsection{Scenario 2} \label{sec:results_scen2}

\begin{figure}[t!]
    \centering
    \includegraphics[width=0.457\textwidth]{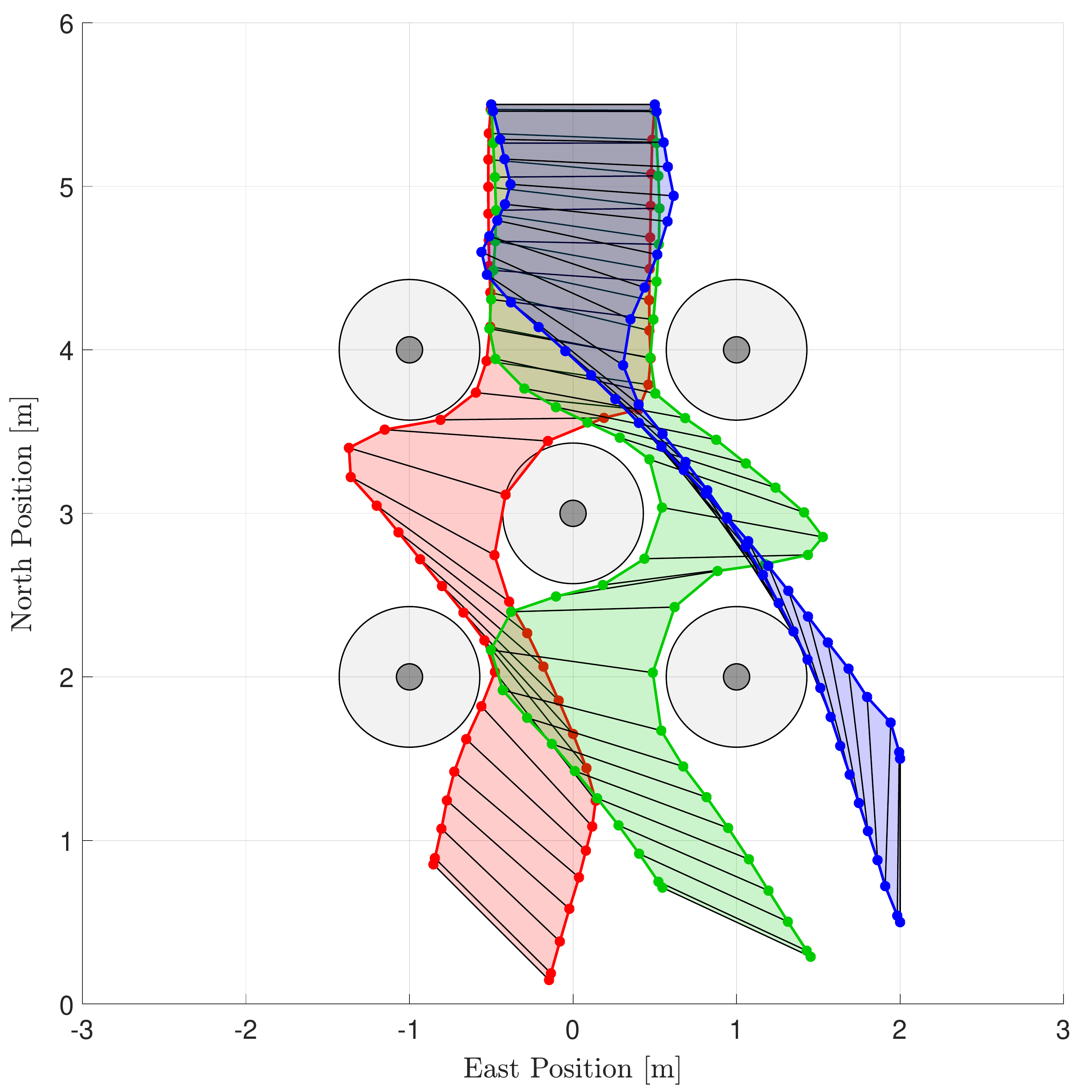}
    \caption{Scenario 2 sample trajectories for $\KK=30$ and a set of $\No=5$ obstacles. Vehicle motion is from top to bottom. Each color represents a separate run, and the dots and black lines represent the vehicles and payload, respectively. The dark gray circles represent the obstacles, while the light gray circle represent the effective keep-out areas. A minor inter-sample constraint violation is observed between the red case and the center obstacle.}
    \vspace{-0.6cm}
    \label{fig:scen2_results}
\end{figure}

For Scenario~2, the performance of Algorithm~\ref{algorithm:scvx} was evaluated by running multiple cases with identical terminal conditions and final times, but with different initial positions, initial formation angles, and obstacle configurations. Example simulations are presented in Figure~\ref{fig:scen2_results}, which shows three sets of trajectories obtained with a temporal resolution of~$\KK=30$ and a manually positioned configuration of~$\No=5$ obstacles. In what follows, recall that Scenario 2 is two-dimensional in nature, and is thus implemented using only two spatial dimensions.

Monte Carlo batch runs were conducted as in~\sref{sec:results_secn1}. The initial positions and angles of the formation were sampled from uniformly random distributions given by~$(r_{i,1}+r_{i,2})/2\in [-2,2]\times [0,1]$ and $\psi\in[-70\degree,70\degree]$. Each case had $\No=4$ obstacles, whose positions $r_{o,\{1,2,3,4\}}$ were sampled from a uniform random distributions of positions~$\roj\in[-2,2]\times[2,4.5]$~m. In generating the obstacle configuration, an inter-obstacle 1-norm distance of at least $\max\{0.8(\Lspacing+2\Wo),2\Wo\}$~m was enforced. This heuristic encouraged the obstacle configuration to be feasibly navigable by the quad-rotor pair. To make sure that the obstacles obstruct a straight path from the initial to the final payload position, one of the obstacles was always placed along the segment from $(r_{i,1}+r_{i,2})/2$ to $(r_{f,1}+r_{f,2})/2$. Other problem parameters were set to the following constant values:
\begin{gather*}
  \tf=4~\text{s},\quad
  r_{f,1}=(-0.5,5.5)~\text{m}, \\
  r_{f,2}=(0.5,5.5)~\text{m},\quad
  \Lspacing=1~\text{m},\quad
  m=0.35~\text{kg}, \\
  \vmax=3~\text{m/s},\quad
  \Tmax=5~\text{N},\quad
  \Ro=0.08~\text{m}, \\
  \Wo = 0.43~\text{m},\quad
  g = 9.81~\text{m/s}^2,\quad
  \Wvc = 10^5\cdot I_{8\times 8}, \\
  \Wtr=50\cdot\blkdiag{I_{8\times 8},0_{4\times 4}}.
\end{gather*}

\begin{table}[t]
    \begin{center}
        \vspace{0.2cm}
    	\caption{Algorithm~\ref{algorithm:scvx} runtime for Scenario 2 [ms].}
 		\begin{tabu} to 0.457\textwidth {c|X[c]X[c]X[c]X[c]X[c]}
 			\hhline{======} \\[-1.50ex]
     		\KK  & Mean & Median & Std. Dev. & Min & Max \\[0.5ex]
     		\hhline{------} \\[-1.50ex]
     		$15$ & $59$ & $53$ & $22$ & $36$ & $157$ \\[1ex]
            $20$ & $138$ & $127$ & $50$ & $67$ & $294$ \\[1ex]
            $25$ & $204$ & $192$ & $76$ & $112$ & $571$ \\[1ex]
            $30$ & $243$ & $212$ & $100$ & $132$ & $693$ \\[1ex]
            \hhline{======}
 		\end{tabu}
        \label{table:scen2_timing}
        \vspace{-0.2cm}
 	\end{center}
\end{table}
\begin{table}[t]
    \begin{center}
    	\caption{Constraint clipping for Scenario 2 [cm].}
 		\begin{tabu} to 0.457\textwidth {c|X[c]X[c]X[c]X[c]X[c]}
 			\hhline{======} \\[-1.50ex]
     		\KK  & Mean & Median & Std. Dev. & Min & Max \\[0.5ex]
     		\hhline{------} \\[-1.50ex]
     		$15$ & $42$ & $42$ & $4$ & $6$ & $43$ \\[1ex]
            $20$ & $5$ & $4$ & $6$ & $0$ & $42$ \\[1ex]
            $25$ & $2$ & $<0.5$ & $2$ & $0$ & $10$ \\[1ex]
            $30$ & $<0.5$ & $0$ & $1$ & $0$ & $5$ \\[1ex]
     		\hhline{======}
 		\end{tabu}
        \label{table:scen2_clipping}
        \vspace{-0.7cm}
 	\end{center}
\end{table}


Tables~\ref{table:scen2_timing} and \ref{table:scen2_clipping} provide runtime and inter-sample constraint violation statistics. For the same reasons given in Scenario~1, runtime increased and constraint clipping decreased for larger $\KK$. However, because Scenario~2 is a larger and more nonlinear problem, Algorithm~\ref{algorithm:scvx} is less capable of running at interactive rates, with a worst-case observed solve rate of approximately 1.5~Hz for $\KK=30$. Furthermore, we note that inter-sample constraint violation occurred almost persistently for $\KK=15$, where the discretized quad-rotor pair was able to ``jump'' through an obstacle in a single discrete time step without violating \eqref{eq:ex2_cstc_cont}. Therefore, in this scenario it appears necessary to either have $\KK\ge 25$ or to decrease $\vmax$ in order to avoid constraint clipping. Finally, 31 infeasible cases were encountered during the Monte-Carlo simulation, all of which can be remedied by increasing $\tf$.




\section{Conclusion \& Future Work} \label{sec:conclusion}

The main contribution of this paper is the application of the recently introduced compound state-triggered constraints to quad-rotor path planning. Two quad-rotor scenarios are outlined, one involving a quad-rotor flying through a hoop, and the second involving two quad-rotors cooperatively negotiating an obstacle course while carrying a beam-like payload. Our simulation results indicate that solutions to both scenarios can be obtained in real-time. Future work will focus on the implementation the proposed methodology on the Autonomous Control Lab's in-house quad-rotor platforms.


\section*{Acknowledgments}

Support for studying the convergence properties of the successive convexification framework was provided by the Office of Naval Research grants N00014-16-1-2877 and N00014-16-1-3144.


\bibliographystyle{IEEEtran}
\bibliography{refs}


\end{document}